\documentclass[final,oneside,letterpaper]{article}%
\usepackage{graphicx}
\usepackage{amsmath}
\usepackage{afterpage}
\usepackage{amstext}
\usepackage{fancyheadings}
\usepackage{latexsym}
\usepackage{plain}
\usepackage{showkeys}
\usepackage{showlabels}
\usepackage{theorem}
\usepackage{wrapfig}
\usepackage{doublespace}
\usepackage{amsfonts}
\usepackage{amssymb}%
\begin{document}

\title{ A BLOWUP PROBLEM OF REACTION DIFFUSION EQUATION RELATED TO THE DIFFUSION
INDUCED BLOWUP\textrm{ PHENOMENON}}
\author{Chu-Pin Lo\\Department of Applied Mathematics,\\Providence University, 200 Chungchi Road, Shalu 43301\\Taichung Hsien, Taiwan. \\cplo@pu.edu.tw}
\maketitle

\begin{abstract}
This work studies nonnegative solutions for the Cauchy, Neumann, and Dirichlet
problems of the logistic type equation $u_{t}=\Delta u+\mu u^{p}-a(x)u^{q}$
with $p,q>1,$ $\mu>0$. The finite time blowup results for nonnegative
solutions under various restrictions on $a(x),p,q,\mu$ are presented. Applying
the results allows one to construct some reaction diffusion systems with the
so-called \textit{diffusion-induced blowup} phenomenon, particularly in the
case of equal diffusion rates. \ \ \ \ \ \ \ \ \ \ \ \ \ \ \ \ \ \ \ \ \ \ \ \ \ \ \ \ \ \ \ \ \ \ \ \ \ \ \ \ \ \ \ \ \ \ \ \ \ \ \ \ \ \ \ 

\end{abstract}

\makeatletter \makeatother \pagestyle{myheadings} \thispagestyle{plain}
\markboth{chu-bin lo}{on the blowup problem of $u_{t}=\Delta
u+\mu u^{p}-a(x)u^{q}$}


\section{Introduction}

\noindent\ \ \ \ \ This study considers the reaction diffusion equation
\begin{equation}
\left\{
\begin{array}
[c]{l}%
u_{t}=\Delta u+\mu u^{p}-a(x)u^{q}\hspace{5mm}\hbox{in}\hspace{3mm}%
D\times(0,T^{\ast}),\\
u(x,0)=u_{0}(x)\geqq0\hspace{5mm}\hbox{on}\hspace{3mm}D,
\end{array}
\right.  \tag{1.1}%
\end{equation}
where $D\subset R^{n}$ is a bounded smooth domain, or $D=R^{n}$, $u_{0}\in
C(\overline{D}),$ and $T^{\ast}$ is the maximum life time of a solution
defined below. When $D$ is a bounded smooth domain, the homogeneous Neumann
boundary condition or Dirichlet boundary condition is imposed:

\bigskip%
\begin{equation}
{\frac{\partial u}{\partial n}}=0,\text{ on }\partial D\times(0,T^{\ast}),
\tag{1.2}%
\end{equation}%
\begin{equation}
{u}=0,\text{ on }\partial D\times(0,T^{\ast}), \tag{1.3}%
\end{equation}
where $n$ denotes the unit exterior normal vector on $\partial D$. It is
assumed that $\mu>0$, $p,q>1$, and $a:D\longrightarrow R$ is a smooth function
which satisfies some additional conditions (see the conditions in the main theorems).

Two motivations exist for studying the above problem.

First, as well known, a classical solution of a reaction diffusion system may
blow up in a finite time (see Levine [\textbf{17}]). Previous literature
usually adopt a spatially homogeneous reaction term. Therefore elucidating the
spatial influence on the blowup problem for a reaction diffusion equation is a
worthwhile task. Fujita [\textbf{8}], Hayakawa [\textbf{12]}, Kobayashi,
Siaro, and Tanaka [\textbf{15}], and Aronson and Weinberger [\textbf{1}]
considered the Cauchy problem for
\[
u_{t}=\Delta u+u^{p},\hspace{5mm}1<p\leqq1+{\frac{2}{n}}.
\]
According to their results, the above problem has no nonnegative nontrivial
solution satisfying
\[
||u(\cdot,t)||_{L^{\infty}}<\infty,\hspace{5mm}\mathit{for}\text{
}\mathit{all}\text{ }t\geqq0.
\]
Thus, the solutions of $u_{t}=\Delta u+\mu u^{p}-a(x)u^{q}$ blow up in finite
time when $a\equiv0.$ On the other hand, if $q>p$ and $a$ is uniformly
positive, the solutions exist all of the time. Therefore, what happens when
$a$ is nonnegative but vanishes at one or more points is worth exploring.

Second, as an application of our blowup result for (1.1) (Theorem 1.2), a
reaction diffusion system, (1.6) below, with the so-called
\textit{diffusion-induced blowup} phenomenon (see [\textbf{20}]) can be
constructed . For a given reaction diffusion system: $u_{t}=D\Delta u+f(x,u),$
with the Cauchy, Neumann or Dirichlet problem, the diffusion term can usually
smooth solutions. However, recent works [\textbf{4,20,21,25,31,32}] indicated
that diffusion terms may cause a blowup in finite time. Such a phenomenon is
called \textit{diffusion-induced blowup. }Notably, for a single equation, no
such phenomenon exists by the comparison principle. Churbanov [\textbf{4}]
provided an example for the Cauchy problem. Meanwhile, Morgan [\textbf{21}],
M. X. Wang [\textbf{31}], Mizoguchi, Ninomiya, and Yanagida [\textbf{20}], and
Weinberger [\textbf{32}] constructed examples for the Neuman and Dirichlet
problems. More precisely, Morgan gave the following example:%

\begin{equation}
\left\{
\begin{array}
[c]{l}%
v_{t}=({\frac{1}{\pi^{2}}})v_{xx}+v,\text{ }0<x<1,\text{ }t>0,\\
u_{t}=u^{2}-|3-v|u^{3},\text{ }0<x<1,\text{ }t>0,\\
v_{x}(0,t)=v_{x}(1,t)=0,\text{ }t>0,\\
v(x,0)=a\mathrm{cos}(\pi x),\hspace{3mm}u(x,0)=u_{0}(x),\text{ }0<x<1,
\end{array}
\right.  \tag{1.4}%
\end{equation}
where $u_{0}$ is a smooth function such that $u_{0}((1/\pi)\arccos(3/a))>0$
with $\left|  a\right|  \geqq3$. Morgan proved that the solutions of (1.4)
blow up in finite time. However, if the diffusion term is excluded, the
solutions of the corresponding ordinary differential equations of (1.4) will
exist globally for all initial data. Later, M. X. Wang [\textbf{31}] modified
the above example by adding the diffusion term in the $u$-equation and letting
$a=3$. In fact, since solution $v$ can be solved explicitly, under appropriate
rescaling of variables $u$ and $t,$ Wang considered the blowup and global
existence problem of the following scalar equation:%

\begin{equation}
\left\{
\begin{array}
[c]{l}%
u_{t}=du_{xx}+u^{p}-(1-\cos\pi x)u^{q},\text{ }0<x<1,\text{ }t>0,\\
u_{x}=0,\text{ }x=0,1;t>0,\\
u(x,0)=u_{0}(x)>0,\text{ }0<x<1,
\end{array}
\right.  \tag{1.5}%
\end{equation}
where $q>p,$ $d>0.$ He proved that if $1<p<q\leqq2p-1$, the solution blows up
in finite time whenever $d$ is small and $u_{0}(x)$ is large. Moreover, the
blowup point occurs at $0$ only. Notably, Wang did not prove the global
existence of the corresponding ordinary equation in this general $p,q$
setting. Therefore, whether or not Wang's example contains the
diffusion-induced blowup phenomenon ($p=2,$ $q=3$ is OK from Morgan's result)
will be verified later. Recently, Mizoguchi \textit{et al. }[\textbf{20}]
provided an example which has the \textit{diffusion-induced blowup} phenomenon
with unequal diffusion rates and asked whether examples with equal diffusion
rates which have such phenomena exist. Weinberger [\textbf{32}] confirmed that
such example exists. This study considers the Cauchy, homogeneous Neumann, and
Dirichlet problems of the following example which has the
\textit{diffusion-induced blowup} phenomenon (see Sec.4):

\textit{%
\begin{equation}
\left\{
\begin{array}
[c]{l}%
v_{t}=d_{1}\Delta v+f(x,v),\mathit{\ \ \ \ \ \ \ \ \ in}\text{ }%
\mathit{D}\times\left(  0,\infty\right)  ,\\
u_{t}=d_{2}\Delta u+\mu u^{p}-|m-v|^{\sigma}u^{q},\mathit{\ }\text{\ \ \ \ }%
\mathit{in}\text{ }\mathit{D}\times\left(  0,\infty\right)  ,\text{ \ }\\
v(x,0)=v_{0}(x)\geqslant0,\hspace{3mm}u(x,0)=u_{0}(x)\geqq0,\hspace{3mm}x\in
D.
\end{array}
\right. \tag{1.6}%
\end{equation}
}where $D=R^{n}$ or a smooth bounded domain, \textit{d}$_{1,}$ $d_{2},$ $\mu$
are nonnegative constants, \ \textit{f}\ is a smooth function, $v(x,t)=v_{0}%
(x)$ independent of $t$ is an entire solution of $R^{n}$ or a solution for the
corresponding homogeneous elliptic Neumann or Dirichlet problem, and
$m=v_{0}(x_{0}),$ for some $x_{0}\in D.$ Notice that since the $v$ equation in
(1.6) does not couple with $u,$ by solving $v$ and substituting it into $u$
equation we can get a scalar equation in the form as in (1.1). Therefore our
blowup results about (1.1) are related to the above reaction diffusion system
(1.6). The corresponding ordinary differential equation is as follows:\textit{%
\begin{equation}
\text{\textit{$\left\{
\begin{array}
[c]{l}%
v_{t}=f(x_{0},v),\\
u_{t}=\mu u^{p}-|m-v|^{\sigma}u^{q},\text{ }t\geqq0,\\
v(0)=\xi\geqq0,u(0)=\eta\geqq0.
\end{array}
\right.  $}}\tag{1.7}%
\end{equation}
}

In sum, our blowup results generalize Wang's results into a more general
equation (in fact, the condition ``$1<p<q\leqq2p-1$'' in Wang's results is the
special case ``$\sigma=2$'' in Theorem 1.2 of this paper) and hold for the
Cauchy, Neumann, and Dirichlet problems. Also, the example herein which
exhibits the \textit{diffusion-induced blowup} phenomenon allows for equal
diffusion rates. Therefore, the question in [\textbf{20}] is also answered positively.

Now, some preparations are made for stating our main results. As well known,
(1.1) with (1.2) or (1.3) in the bounded smooth domain $D$ has an unique local
classical solution $u(x,t)$ defined in some maximal interval of existence
$0<t<T^{\ast}(u_{0})$ such that
\[
\lim_{t\rightarrow T^{\ast}(u_{0})}||u(\cdot,t)||_{\infty}=\infty
\]
whenever $T^{\ast}(u_{0})<\infty$ (see Henry [\textbf{13]}). $T^{\ast}(u_{0})$
is called the blowup time if $T^{\ast}(u_{0})<\infty$. If there exist
sequences $\{x_{k}\}$, $\{t_{k}\}$ with $t_{k}\in(0,T^{\ast})$, $x_{k}%
\in\overline{D}$ satisfying
\[
\lim_{k\rightarrow\infty}x_{k}=x_{0},\text{ }\lim_{k\rightarrow\infty}%
t_{k}=T^{\ast},\hspace{3mm}\lim_{k\rightarrow\infty}|u(x_{k},t_{k})|=\infty,
\]
this $x_{0}$ is called a blowup point. As for the Cauchy problem (i.e.,
$D=R^{n})$, when $u_{0}$ is bounded uniformly continuous, a local solution in
the class of bounded uniformly continuous functions is guaranteed (see
[\textbf{18}, Sect. 1.4.]), and the blowup time can be similarly defined as in
the bounded smooth domain case. This work considers solutions for the Cauchy
problem in the more general space:%

\begin{equation}
|u(x,t)|\leqq M|x|^{\frac{2}{p-1}},\hspace{3mm}(x,t)\in R^{n}\times(0,T^{\ast
}),\hspace{3mm}M>0. \tag{1.8}%
\end{equation}

Where the Cauchy problem for (1.1) is denoted by $(\mathcal{I)}$, and the
homogeneous Neumann and Dirichlet problem is represented by $(\mathcal{II)}$
and ($\mathcal{III}\mathbb{)}$ respectively. The solution of $(\mathcal{I)}$
is always assumed to satisfy (1.8).

Before Proposition 1.1, i.e., the key lemma in this study, is stated, some
heuristic reasons for this proposition are explained. Notably, the equation in
(1.1) with $a(x)=|x|^{\sigma}$ is invariant under the scaling: $t\rightarrow
\lambda t$, $x\rightarrow\lambda^{\frac{1}{2}}x$, $u\rightarrow\lambda
^{\frac{-1}{p-1}}u$, if $\sigma=2(q-p)(p-1)^{-1}$ (or $q=(\frac{\sigma}%
{2}+1)p-\frac{\sigma}{2}$). This feature suggests finding the self-similar
lower solution and using it as a blowing up lower solution to obtain the
blowup results. The following problem (1.9) is obtained by using similarity
variables (see Lemma 3.1).

\bigskip

\bigskip

\textbf{Proposition 1.1. }\textit{Let $q\geqq p>1$ and $\sigma\geqq0$. When
$n\geqq3,$ further assume that $p<{\frac{n+2}{n-2}}$}. \textit{Then, there
exists $\mu_{0}>0$ such that the Dirichlet problem
\begin{equation}
\left\{
\begin{array}
[c]{l}%
\Delta w(y)-{\frac{1}{2}}y\cdot\nabla w(y)-{\frac{1}{p-1}}w(y)+\mu_{0}%
w^{p}(y)-|y|^{\sigma}w^{q}(y)=0,\\
w|_{\partial B(0,r_{0})}=0,
\end{array}
\right. \tag{1.9}%
\end{equation}
has a radially symmetry positive weak solution $w_{0}\in H_{0}^{1}%
(B(0,r_{0}),\rho)$, where }$\rho(y)$$=$$exp({\frac{-|y|^{2}}{4}}),$
$B(0,r_{0})$ \textit{is the} \textit{open} \textit{ball centered at }%
$0$\textit{ with radius} $r_{0},$ \textit{and} $H_{0}^{1}(B(0,r_{0}),\rho)$
\textit{is the} $H_{0}^{1}$ \textit{Sobolev space with weight }$\rho(y)$
(\textit{see Escobedo and Kavian }[\textbf{5}]).\textit{ Moreover, if
$1<q<{\frac{n+2}{n-2}}$ when $n\geqq3$, then $w_{0}\in C^{2+\alpha}%
(B(0,r_{0}))$ for some $0<\alpha<1$. }\newline

\bigskip

\bigskip

Proposition 1.1 can be proved similarly to Proposition 2 of [\textbf{9}]
(variational method) and Theorem 3.12 of [\textbf{5}] (using the bootstrap
method for the regularity). Therefore the proof of Proposition 1.1 is omitted here.

Define
\begin{equation}
w(r)\equiv\left\{
\begin{array}
[c]{l}%
w_{0}(r),\hspace{5mm}0\leqq r\leqq r_{0},\\
0,\hspace{5mm}r>r_{0}.
\end{array}
\right.  \tag{1.10}%
\end{equation}
The following theorem is the main result concerning blowup:

\bigskip

\bigskip

\textbf{Theorem 1.2. }\textit{Assume that $u(x,t)$ is a nonnegative classical
upper solution of ($\mathcal{I)}$}.\textit{ Also assume $q\geqq p>1$ and
$a(x)$ is a real valued smooth function satisfying
\[
a(x)\leqq M|x-x_{0}|^{\sigma},\hspace{3mm}x\in R^{n},\hspace{3mm}%
M>0,\hspace{3mm}\sigma\geqq2(q-p)(p-1)^{-1},
\]
where $a(x_{0})=0$. When $n\geqq3,$ it is further assumed that }$q<{\frac
{n+2}{n-2}.}$ \textit{Let $\mu\geqq M^{\frac{p-1}{q-1}}\mu_{0}$ in $(1.1)$ and
$t_{0}$ be any negative number where $\sigma=2(q-p)(p-1)^{-1}$ or }$-1\leqq
t_{0}<0$\textit{\ where $\sigma>2(q-p)(p-1)^{-1}$. If $u_{0}(x)\geqq
M^{-{\frac{1}{q-1}}}(-t_{0})^{\frac{-1}{p-1}}w({\frac{|x-x_{0}|}{\sqrt{-t_{0}%
}}})$ for all $x\in R^{n}$, then $u(x,t)$ will blow up at a finite time which
is before or equal to ``$-t_{0}$''. Similarly, if $u(x,t)$ is a nonnegative
classical upper solution of ($\mathcal{II)}$, or ($\mathcal{III})$ in $D\equiv
B(x_{0},r_{0}\sqrt{-t_{0}})$ with }$r_{0\text{ \ }}$\textit{as in Proposition
1.1, then the above result also holds under the same conditions as in
($\mathcal{I)}$.}

\bigskip

\bigskip

\textbf{Remark.} (i) The proof of Theorem 1.2 indicates that the lower bound
of blowup rate is $(T^{\ast}-t)^{\frac{-1}{p-1}}$. If $a(x)$ is a nonnegative
term and $u(x,t)$ is a nonnegative classical solution of the Cauchy problem
$(\mathcal{I)}$, then the upper bound of the blowup rate is also $(T^{\ast
}-t)^{\frac{-1}{p-1}}$ by [\textbf{10}, Theorem 3.7] and the comparison
principle. Therefore, the blowup rate is $(T^{\ast}-t)^{\frac{-1}{p-1}}$.

(ii) The result obtained in Theorem 1.2 is natural from a biological
perspective. Since the term, [$\mu a(x)^{-1}]^{\frac{1}{q-p}},$ is explained
as the saturation level of the species, the conditions of $\mu,$ $a(x)$ in
Theorem 1.2 imply [$\mu a(x)^{-1}]^{\frac{1}{q-p}}$ will become very large,
particularly at the zero points of \ $a(x)$, and the population density tends
to increase rapidly at such points, as stated in Theorem 1.2.

\bigskip

\bigskip

The structure of the blowup set is as follows:

\bigskip

\bigskip

\textbf{Theorem 1.3. }\textit{Assume that}\textbf{ }$q>p>1$ \textit{and}
\textit{a(x)} \textit{is a nonnegative smooth function. For problem}
($\mathcal{III}$), \textit{if }$\sigma\geqq2(q-p)(p-1)^{-1},$ \textit{where}
$\sigma$\textit{ is the zero order of }$a(x)$ \textit{at any zero point, then
the blowup points} \textit{occur at the zero points of a(x) only. For problem}
($\mathcal{II}$), \textit{if a(x) satisfies the following condition: for every
x}$\in\partial D$, \textit{either }$\frac{\partial a(x)}{\partial n}\leqq0,$
\textit{or }$a(x)=0,$ \textit{then the blowup points occur at the zero points
of a(x) only}.

\bigskip

\bigskip

The methods applied herein to verify the accuracy of these results are the
similarity variables method and comparison principle in the weak or classical
sense. The similarity variables method is widely applied to study the blow-up
or global existence phenomenon (see [\textbf{5, 9, 10, 24, 28}]). Generally,
two kinds of self-similar solutions exist, backward and forward. The backward
self-similar solution is relevant to the blowup problem, while the forward
self-similar solution is relevant to the global existence problem. Proposition
1.1 provides a backward self-similar classical solution of the original Eqn.
(1.1) with $a(x)=|x|^{\sigma}$ in some special domain. To apply a comparison
principle in the weak sense for the Cauchy problem (Lemma 2.3) or the Neumann
and Dirichlet problems (Lemma 2.2) to the subject domain $D$, in Lemma 3.2 the
above backward self-similar solution is suitably extended to the whole domain
$D$ via Lemma 2.4, which is a connecting lemma for producing weak lower
solution in a larger domain. Thus the blowup result (Theorem 1.2) follows from
comparison with the weak lower solution in $D$.

The rest of this paper is organized as follows. Section 2 lists some
comparison principles in a weak sense. Section 3 proves the blowup results,
Theorem 1.2 and Theorem 1.3. Finally, Section 4 presents some examples that
illustrate the \textit{diffusion-induced blowup} phenomenon. \ \ \ \ \ \ \ \ \ \ \ 

\bigskip

\bigskip

\section{The Comparison Principle}

The main tools applied within for proving the main results are the comparison
principles in a weak sense (see [\textbf{7, 16, 19,} \textbf{30}]). In this
section, these comparison principles are stated or derived in the suitable
form convenient for the application to the Dirichlet, Neumann and Cauchy
problems herein.

Let $\Omega\subset R^{n+1}$ be a bounded smooth domain. The following linear
equation is considered.%

\begin{equation}
{\frac{\partial u}{\partial t}}\equiv\Delta u+\sum_{i=1}^{n}a_{i}%
(x,t){\frac{\partial u}{\partial x_{i}}}+a(x,t)u\hspace{5mm}\mbox{in}\hspace
{3mm}\Omega, \tag{2.1}%
\end{equation}
where $a_{j}$, ${\frac{\partial a_{i}}{\partial x_{k}}}$, $a$ are
H$\mathrm{\ddot{o}}$lder continuous in $(x,t)\in\Omega$ and $a\leqq0$ in
$\Omega$. The so-called Strong Maximum Principle in a weak sense is as follows
(see [\textbf{7}]).

\bigskip

\bigskip

\textbf{Lemma 2.1. }\textit{Assume that $u(x,t)$ is a bounded measurable
function in $\Omega$ satisfying
\begin{equation}
\int_{E}u(x,t)(\Delta v-\sum_{i}{\frac{\partial(a_{i}v)}{\partial x_{i}}%
}+{\frac{\partial v}{\partial t}}+av)dxdt\geqq0,\tag{2.2}%
\end{equation}
for every compact subset $E\subset\Omega$ and nonnegative test function $v$
with $supp(v)$$\subset E$. Assume that $u(x^{0},t^{0})=$ $\underset
{(x,t)\in\Omega}{\mathit{\mathit{esssup}}}$ $u(x,t)\geqq0$ with $(x^{0}%
,t^{0})\in\Omega$. Use $C(x^{0},t^{0})$ to denote the set of all point
$(x^{1},t^{1})$ in $\Omega$ that can be connected to $(x^{0},t^{0})$ by a
smooth curve in $\Omega$ along which the $t$-coordinate is increasing. If
$u|_{C(x^{0},t^{0})}$ is a continuous function, then $u\equiv M$ almost
everywhere in $C(x^{0},t^{0})$.}

\bigskip

\bigskip

Lemma 2.1 is used herein to prove a comparison principle for the Neumann and
Dirichlet problems. Let $\Omega\equiv B(0,r_{0}\sqrt{-t_{0}})\times
(t_{0},t_{1})$ with $t_{0}<t_{1}<0$ fixed, $S\equiv\partial B(0,r_{0}%
\sqrt{-t_{0}})\times(t_{0},t_{1}]$, and $I\equiv\{(x,t)\in\Omega
:|x|>r_{0}\sqrt{-t}\}$.

\bigskip

\bigskip

\textbf{Lemma 2.2. }\textit{Assume that}\textbf{ }\textit{\ $u(x,t)$ is a
continuous function in $\overline{\Omega}$. Let $h(x,t)\in C^{\alpha
,{\frac{\alpha}{2}}}(\Omega)$ be a bounded function. Assume that $-u$
satisfies (2.2) with }$a_{i}$ \textit{and} $a$ \textit{replaced by} $0$
\textit{and} $-h$ \textit{respectively.} \textit{For} \textit{the Dirichlet}
\textit{problem}, \ \textit{if ${u}|_{S}\geqq0$ and $u(x,t_{0})\geqq0$, then
$u(x,t)\geqq0$ in $\Omega$. For the Neumann problem, further assume that
$u|_{I}$ is smooth. Then under }$\frac{\partial u}{\partial n}|_{S}\geqq0$
\textit{and }$u$\textit{$(x,t_{0})\geqq0,$ solution $u(x,t)$ is also
nonnegative in }$\Omega.$

\bigskip

\bigskip

\textit{Proof.\hspace{2mm}} Since $-u$ satisfies (2.2), we have
\[
\int_{\Omega}u(-\Delta v-v_{t}+h(x,t)v)\geqq0
\]
for all nonnegative test function with $supp(v)$$\subset\subset$$\Omega$. Let
$\nu\equiv ve^{lt}$ with $l>-h$. Then
\[
-\Delta v-v_{t}+h(x,t)v=[-\Delta\nu+(h+l)\nu-\nu_{t}]e^{-lt}.
\]
Hence
\[
\int_{\Omega}ue^{-lt}[-\Delta\nu-\nu_{t}+(h+l)\nu]dxdt\geqq0.
\]
Set $u^{\ast}(x,t)\equiv u(x,t)e^{-lt}$. To prove the result, it suffices to
show $u^{\ast}(x,t)\geqq0$. By $u(x,t_{0})\geqq0$, $u^{\ast}(x,t_{0})\geqq0$
is obtained. Assume that the minimum of $u^{\ast}$ over $\overline{\Omega}$ is
negative, say $m$, and $u^{\ast}(\overline{x},\overline{t})=m$ for some point
$(\overline{x},\overline{t})\in\overline{\Omega}$. First, if $(\overline
{x},\overline{t})$$\in$$\Omega\cup$$B(0,r_{0}\sqrt{-t_{0}})\times\{t_{1}\}$,
then Lemma 2.1 can be used to obtain $u^{\ast}(x,t_{0})\equiv m$. This
contradicts the fact ``$u^{\ast}(x,t_{0})\geqq0$.'' Similarly, $\overline{t}$
is not equal to $t_{0}$. For the Dirichlet problem, $(\overline{x}%
,\overline{t})$ does not occur at $S,$ producing a contradiction. Therefore we
get a contradiction. For the Neumann problem, if $(\overline{x},\overline
{t})\in S$, then since $\min_{\overline{I}}u^{\ast}=m$, $u^{\ast}|_{I}$ is
smooth, and $u^{\ast}(x,t)|_{I}>m$ (by the previous argument), ${\frac
{\partial u^{\ast}}{\partial n}}(\overline{x},\overline{t})<0$ on $S$ \ is
obtained by the classical Hopf maximum principle. This finding also
contradicts ``${\frac{\partial u}{\partial n}}|_{S}\geqq0$.'' The proof is
completed.\hfill\textbf{$\blacksquare$}

\bigskip

\bigskip

Consider the following problem
\begin{equation}
\left\{
\begin{array}
[c]{l}%
u_{t}=\Delta u+f(x,t,u)\hspace{5mm}\mbox{in}\hspace{3mm}R^{n}\times(0,T^{\ast
}),\\
u|_{t=0}=u_{0}(x)\hspace{5mm}\mbox{in}\hspace{3mm}R^{n},
\end{array}
\right.  \tag{2.3}%
\end{equation}
where $T^{\ast}>0$ and $u_{0}(x)$ is a continuous function. For the above
Cauchy problem, the corresponding comparison principle is deduced in
[\textbf{30}]:

\bigskip

\bigskip

\textbf{Lemma 2.3. }\textit{Assume that $\overline{u}(x,t)$($\underline
{u}(x,t)$) is a continuous function on $R^{n}\times\lbrack0,T^{\ast})$
satisfying
\[
\overline{u}|_{t=0}\geqq(\leqq)u_{0}(x)
\]
and
\begin{align}
&  \int_{R^{n}}\overline{u}(x,t)\eta(x,t)dx|_{t=0}^{t=T_{1}}\tag{2.4}\\
&  \geqq(\leqq)\int_{0}^{T_{1}}\int_{R^{n}}[\overline{u}(x,s)(\Delta\eta
+\eta_{t})(x,s)+\eta(x,s)f(x,s,\overline{u})]dxds\nonumber\\
&  \hspace{5mm}\mbox{for}\hspace{3mm}T_{1}\in\lbrack0,T^{\ast}),\nonumber
\end{align}
where $\eta$ is a nonnegative, smooth function on $R^{n}\times\lbrack
0,T^{\ast})$ with $supp(\eta(\cdot,t))$$\subset\subset$$R^{n}$ for all
$t\in\lbrack0,T^{\ast})$. Assume $(\overline{u}-\underline{u})(x,t)\geqq
-Bexp(\beta|x|^{2})$ on $R^{n}\times(0,T^{\ast})$ and
\[
f(x,t,\overline{u}(x,t))-f(x,t,\underline{u}(x,t))\geqq c(x,t)(\overline
{u}-\underline{u})(x,t),
\]
where $B$, $\beta$$>0$, $c\in$$C_{loc}^{\alpha,{\frac{\alpha}{2}}}(R^{n}%
\times(0,T^{\ast}))$, and $c(x,t)\leqq c_{0}(|x|^{2}+1)$ on $R^{n}%
\times(0,T^{\ast})$ for some $c_{0}$. Then, $\overline{u}\geqq\underline{u}$
on $R^{n}\times\lbrack0,T^{\ast})$.}

\bigskip

\bigskip

\textbf{Remark. }The $-u$ ($u$) in Lemma 2.2 is called a continuous weak lower
(upper) solution of the Dirichlet problem, or the Neumann problem. Also,
$\overline{u}$ and $\underline{u}$ in Lemma 2.3 will be called the continuous
weak upper and lower solutions of the Cauchy problem ($2.3$) respectively.

\bigskip

\bigskip

For later convenience, we now prove the following ``\textit{connecting
lemma}'' which is used to combine two classical lower solutions into one
continuous weak lower solution. Let $\Omega_{1}$ be a bounded domain in
$R^{n}\times(t_{0},t_{1})$, $\Omega_{2}$$\equiv$$R^{n}\times(t_{0}%
,t_{1})-\overline{\Omega}_{1}$, $D_{i}(t)$$\equiv$ $\{x\in R^{n}%
:(x,t)\in\Omega_{i}\}$ for $t\in(t_{0},t_{1})$, and $D\equiv$$\bigcup
_{t\in(t_{0},t_{1})}D_{1}(t)$.

\bigskip

\bigskip

\textbf{Lemma 2.4. }\textit{Let $\underline{u}(x,t)$ be a smooth function and
satisfy
\begin{equation}
u_{t}\leqq\Delta u+f(x,t,u)\tag{2.5}%
\end{equation}
on $\Omega_{1}$, $\Omega_{2}$ respectively. If $\underline{u}$ is continuous
on $R^{n}\times\lbrack t_{0},t_{1}]$ and ${\frac{\partial\underline{u}%
(\cdot,t)}{\partial n_{1}}}|_{\partial D_{1}(t)}$$+$${\frac{\partial
\underline{u}(\cdot,t)}{\partial n_{2}}}|_{\partial D_{2}(t)}$ $\leqq0$, then
$\underline{u}$ satisfies (2.4) in the case ``}$\leqq$\textit{'' with (}%
$0$\textit{,}$T^{\ast}$\textit{) replaced by }$(t_{0,}t_{1})$\textit{, where
$n_{i}$ is the unit outer-normal vector of $D_{i}(t)$ for all $t\in
(t_{0},t_{1})$.}

\bigskip

\bigskip

\textit{Proof.\hspace{2mm}} For a nonnegative, smooth function $\phi$ on
$R^{n}\times\lbrack t_{0},t_{1}]$ with $supp(\phi(\cdot,t))$$\subset\subset
$$R^{n}$ for all $t\in\lbrack t_{0},t_{1}]$, $(\Delta)$ denotes the following
term
\begin{align*}
(\Delta)  &  \equiv\int_{t_{0}}^{t_{1}}\int_{R^{n}}[\underline{u}_{t}%
-\Delta\underline{u}-f(x,t,\underline{u})]\phi\\
&  =\int_{\Omega_{1}}[\underline{u}_{t}-\Delta\underline{u}-f(x,t,\underline
{u})]\phi+\int_{\Omega_{2}}[\underline{u}_{t}-\Delta\underline{u}%
-f(x,t,\underline{u})]\phi.
\end{align*}
Let $(A)$$\equiv$$\int_{\Omega_{1}}{\underline{u}_{t}}\phi$, $(B)$$\equiv
$$\int_{\Omega_{2}}{\underline{u}_{t}}\phi$, $(C)$$\equiv$$\int_{\Omega_{1}%
}{\Delta\underline{u}}\phi$, and $(D)$$\equiv$$\int_{\Omega_{2}}%
{\Delta\underline{u}}\phi$.
\[
(A)=\int_{D}dx\int_{a(x)}^{b(x)}{\underline{u}_{t}}\phi dt=\int_{D}%
\underline{u}\phi|_{a(x)}^{b(x)}dx-\int_{\Omega_{1}}\underline{u}\phi_{t},
\]%
\begin{align*}
(B)  &  =\int_{D}dx\int_{[t_{0},t_{1}]-[a(x),b(x)]}{\underline{u}_{t}}\phi
dt+\int_{R^{n}-D}dx\int_{t_{0}}^{t_{1}}{\underline{u}_{t}}\phi dt\\
&  =\int_{D}(\underline{u}\phi|_{b(x)}^{t_{1}}+\underline{u}\phi|_{t_{0}%
}^{a(x)})dx+\int_{R^{n}-D}\underline{u}\phi|_{t_{0}}^{t_{1}}dx-\int
_{\Omega_{2}}\underline{u}\phi_{t}.
\end{align*}
Adding $(A)$ and $(B)$ leads to
\[
(A)+(B)=-\int_{R^{n}\times\lbrack t_{0},t_{1}]}\underline{u}\phi_{t}%
+\int_{R^{n}}\underline{u}\phi|_{t_{0}}^{t_{1}}dx.
\]%
\begin{align*}
(C)  &  =\int_{t_{0}}^{t_{1}}dt\int_{D_{1}(t)}\Delta\underline{u}\phi dx\\
&  =\int_{t_{0}}^{t_{1}}dt\int_{\partial D_{1}(t)}\phi{\frac{\partial
\underline{u}}{\partial n_{1}}}d\sigma-\int_{t_{0}}^{t_{1}}dt\int_{\partial
D_{1}(t)}\underline{u}{\frac{\partial\phi}{\partial n_{1}}}d\sigma+\int
_{t_{0}}^{t_{1}}dt\int_{D_{1}(t)}\underline{u}\Delta\phi dx,
\end{align*}%
\begin{align*}
(D)  &  =\int_{t_{0}}^{t_{1}}dt\int_{D_{2}(t)}\Delta\underline{u}\phi dx\\
&  =\int_{t_{0}}^{t_{1}}dt\int_{\partial D_{2}(t)}\phi{\frac{\partial
\underline{u}}{\partial n_{2}}}d\sigma-\int_{t_{0}}^{t_{1}}dt\int_{\partial
D_{2}(t)}\underline{u}{\frac{\partial\phi}{\partial n_{2}}}d\sigma+\int
_{t_{0}}^{t_{1}}dt\int_{D_{2}(t)}\underline{u}\Delta\phi dx.
\end{align*}
Substituting the above results into $(\Delta)$ produces
\begin{align*}
(\Delta)  &  =\int_{R^{n}}\underline{u}\phi|_{t_{0}}^{t_{1}}dx-\int
_{R^{n}\times\lbrack t_{0},t_{1}]}\underline{u}\phi_{t}+\underline{u}%
\Delta\phi+f(x,t,u)\phi\\
&  +\int_{t_{0}}^{t_{1}}dt(\int_{\partial D_{1}(t)}\underline{u}%
{\frac{\partial\phi}{\partial n_{1}}}d\sigma+\int_{\partial D_{2}%
(t)}\underline{u}{\frac{\partial\phi}{\partial n_{2}}}d\sigma)\\
&  -\int_{t_{0}}^{t_{1}}dt(\int_{\partial D_{1}(t)}\phi{\frac{\partial
\underline{u}}{\partial n_{1}}}d\sigma+\int_{\partial D_{2}(t)}{\frac
{\partial\underline{u}}{\partial n_{2}}}\phi d\sigma).
\end{align*}
Since $\underline{u}(x,t)$ is a classical lower solution on $\Omega_{1}$,
$\Omega_{2},$ respectively and $\phi\geqq0$, $(\Delta)\leqq0$ is obtained. By
the smoothness of $\phi$ and ${\frac{\partial\underline{u}(\cdot,t)}{\partial
n_{1}}}|_{\partial D_{1}(t)}$$+$${\frac{\partial\underline{u}(\cdot
,t)}{\partial n_{2}}}|_{\partial D_{2}(t)}\leqq0$, the proof is
completed.\hfill\textbf{$\blacksquare$}

\bigskip

\bigskip

\textbf{Remark. }Lemma 2.4 can also be proved by modifying the argument of the
proof of Lemma 3.3 in [\textbf{19}]. A slightly different proof is presented
herein for completeness.

\bigskip

\bigskip

\section{Proofs of Blowup Results}

In this section, the main results about the blowup problem, Theorem 1.2 and
Theorem 1.3, will be proved. Since the equation in (1.1) with
$a(x)=|x|^{\sigma}$ has some scaling invariant property under the condition
$q=({\frac{\sigma}{2}}+1)p-{\frac{\sigma}{2}}$, the following lemma states the
relation between the original equation and the equation which is satisfied by
the self-similar lower or upper solution.

\bigskip

\bigskip

\textbf{Lemma 3.1. }\textit{Assume that $\mu>0$, $\sigma\geqq0$, $p>1$, and
$q=({\frac{\sigma}{2}}+1)p-{\frac{\sigma}{2}}+m$, for some $m\in R$. Let
$\Omega_{1}$ be a domain of $R^{n}\times(t_{1},0)$ for some $t_{1}<0$ and
$Q_{1}$ be defined by $Q_{1}\equiv\{{\frac{|x|}{\sqrt{-t}}}:(x,t)\in\Omega
_{1}\}$. Assume that $w(r)$ is a smooth function defined on $Q_{1}$ and
$u(x,t)$ is defined by $u(x,t)\equiv(-t)^{\frac{-1}{p-1}}w({\frac{|x|}%
{\sqrt{-t}}})$ on $\Omega_{1}$. The following is then obtained:
\begin{equation}
u_{t}\leqq\Delta u+\mu u^{p}-|x|^{\sigma}u^{q}\hspace{5mm}on\hspace{3mm}%
\Omega_{1}\tag{3.1}%
\end{equation}
if and only if
\begin{equation}%
\begin{array}
[c]{l}%
w_{rr}+{\frac{n-1}{r}}w_{r}-{\frac{1}{2}}rw_{r}-{\frac{1}{p-1}}w+\mu
w^{p}-r^{\sigma}(-t)^{\frac{-m}{p-1}}w^{q}\geqq0\\
\hspace{5mm}{on}\hspace{3mm}Q_{1},\hspace{3mm}where\hspace{3mm}r={\frac
{|x|}{\sqrt{-t}}}.
\end{array}
\tag{3.2}%
\end{equation}
}

\bigskip

\bigskip

\textit{Proof.\hspace{2mm}} The proof is straightforward.\hfill
\textbf{$\blacksquare$}

\bigskip

\bigskip

From the above lemma, the backward self-similar solution satisfies the
equation in $(1.9)$ or $(3.2)$ when the solution is radially symmetric.

To apply $w_{0}$ in Proposition 1.1 to our domain, $w_{0}$ is extended by
(1.10) and then the following lemma is obtained:

\bigskip

\bigskip

\textbf{Lemma 3.2. }\textit{Suppose $q\geqq p>1$ and $\sigma\geqq
2(q-p)(p-1)^{-1}$. When n}$\geqq3,$ \textit{it is further assumed that}
\textit{$1<q<{\frac{n+2}{n-2}.}$ Let $w(r)$ be the function defined by
$(1.10)$ and then $\underline{u}$ is defined as follows:
\[
\underline{u}(x,t)\equiv(-t)^{\frac{-1}{p-1}}w({\frac{|x|}{\sqrt{-t}}}%
)\hspace{5mm}\mbox{for}\hspace{5mm}(x,t)\in R^{n}\times\lbrack t_{0},t_{1}],
\]
where $t_{0}<t_{1}<0$, if $\sigma=2(q-p)(p-1)^{-1}$ or -1}$\leqq$
\textit{$t_{0}<t_{1}<0$, if $\sigma>2(q-p)(p-1)^{-1}$. If $u_{0}%
(x)\geqq\underline{u}(x,t_{0})$ and $\mu\geqq\mu_{0}$, then $\underline
{u}(x,t)$ is a continuous weak lower solution of }$\ problems\mathit{\ }%
$($\mathcal{I}$), ($\mathcal{II}$), \textit{and} ($\mathcal{III}$)
\textit{with} ($0,T^{\ast}$), $a(x),$ \textit{and }$D$ \textit{replaced by}
($t_{0}$, $t_{1}$), $|x|^{\sigma}$, and $B(0,r_{0}\sqrt{-t_{0}}),$
\textit{respectively, where }$r_{0}$ \textit{is defined in}
\textit{Proposition 1.1}. \textit{ }

\bigskip

\bigskip

\textit{Proof.\hspace{2mm}} As mensioned in Proposition \textit{1.1,}
$w_{0}(y)$ is a smooth, nonnegative, radially symmetric, and nonconstant
solution of (1.9). By the maximal principle, $w_{0}(y)$ is positive on
$B(0,r_{0})$ and ${\frac{\partial w_{0}}{\partial n}}|_{\partial B(0,r_{0}%
)}={w_{0}}^{^{\prime}}(r_{0})<0$. Thus, our result follows from (1.10), Lemma
2.4, and Lemma 3.1 whenever $\mu\geqq\mu_{0}$. \hfill\textbf{$\blacksquare$}

\bigskip

\bigskip

Applying the comparison principle for the Cauchy problem (Lemma 2.3) to
problem ($\mathcal{I}$) allows us to obtain the following lemma:

\bigskip

\bigskip

\textbf{Lemma 3.3.}

\begin{itemize}
\item[(a)] \textit{\ Assume that $\overline{u}(x,t)$ is a continuous weak
upper solution of ($\mathcal{I}$). If $\overline{u}(x,t)\in C_{loc}%
^{\alpha,{\frac{\alpha}{2}}}(R^{n}\times(t_{0},t_{1}))$, then $\overline
{u}(x,t)\geqq0$ on $R^{n}\times\lbrack t_{0},t_{1}]$.}

\item[(b)] \textit{If\hspace{3mm}$u(x,t)$ is a classical upper solution of
($\mathcal{I}$) and $\underline{u}(x,t)$ is the continuous weak lower solution
in Lemma 3.2 , then $u(x,t)\geqq\underline{u}(x,t)$ on $R^{n}\times\lbrack
t_{0},t_{1}]$ .}
\end{itemize}

\bigskip

\bigskip

Similarly, applying the comparison principle for the Neumann and Dirichlet
problems (Lemma 2.2) to ($\mathcal{II}$),($\mathcal{III}$) allows us to obtain
the following lemma:

\bigskip

\bigskip

\textbf{Lemma 3.4. }\textit{Let $u$ be a nonnegative, bounded classical upper
solution of ($\mathcal{II}$) or ($\mathcal{III}$) and $\underline{u}(x,t)$ be
defined in Lemma 3.2. Then, $u(x,t)\geqq\underline{u}(x,t)$ in $B(0,r_{0}%
\sqrt{-t_{0}})\times(t_{0},t_{1})$. }

\bigskip

\bigskip

\bigskip

Now, the special case for Theorem 1.2 is proved, i.e., $a(x)=|x|^{\sigma}$ in
$(1.1)$.

\bigskip

\bigskip

\textbf{Lemma 3.5. }\textit{Assume that $u(x,t)$ is a nonnegative classical
upper solution of $(\mathcal{I)}$. Suppose $q\geqq p>1$ and $a(x)=|x|^{\sigma
}$ with $\sigma\geqq$ }$2(q-p)(p-1)^{-1}$\textit{. When n}$\geqq3,$ \textit{it
is further assumed that} \textit{$1<q<{\frac{n+2}{n-2}.}$ Let $\mu\geqq\mu
_{0}$ in $(1.1)$ and $t_{0}$ be any negative number when $\sigma
=2(q-p)(p-1)^{-1}$ or -1}$\leqq$ \textit{$t_{0}<0$ when $\sigma
>2(q-p)(p-1)^{-1}$. If $u_{0}(x)\geqq(-t_{0})^{\frac{-1}{p-1}}w({\frac
{|x|}{\sqrt{-t_{0}}}})$, then $u(x,t)$ will blow up at finite time which is
before or equal to ``$-t_{0}$''.\newline Similarly, if $u(x,t)$ is a
nonnegative classical upper solution of $(\mathcal{II)}$ or $(\mathcal{III)}$
in $D=B(0,r_{0}\sqrt{-t_{0}})$, then the above result also holds under the
same conditions. }\newline

\bigskip

\bigskip

\textit{Proof.\hspace{2mm}} First, problem $(\mathcal{I)}$ is considered. By
translation, the result of Lemma 3.3 also holds for nonnegative time, i.e.,
$u(x,t)\geqq\underline{u}(x,t+t_{0})$ for $(x,t)\in R^{n}\times\lbrack
0,-t_{0})$. Since
\[
\underline{u}(x,t+t_{0})=(-(t+t_{0}))^{\frac{-1}{p-1}}w({\frac{|x|}%
{\sqrt{-(t+t_{0})}}}),
\]%
\[
\underline{u}(0,t+t_{0})=(-(t+t_{0}))^{\frac{-1}{p-1}}w_{0}(0),\hspace
{3mm}\mbox{and}
\]%
\[
\lim_{t\rightarrow(-t_{0})^{-}}\underline{u}(0,t+t_{0})=\infty,
\]
$u(x,t)$ will blow up at some finite time which is before or equal to
``$-t_{0}$''.

For problems $(\mathcal{II)}$, $(\mathcal{III)},$ the results can be proved by
Lemma 3.4 similarly to the above proof.\hfill\textbf{$\blacksquare$}

\bigskip

\bigskip

\textbf{Proof of Theorem 1.2. }By rescaling: $u^{\ast}(x,t)\equiv M^{\frac
{1}{q-1}}u(x,t)$ and applying Lemma 3.5 to $u^{\ast}$, the theorem can easily
be proved.\hfill\textbf{$\blacksquare$}

\bigskip

\bigskip

\bigskip\textbf{Proof of Theorem 1.3. }First the problem\textbf{
(}$\mathcal{II}$\textbf{) }is considered. Without loss of generality, $n=1$ is
assumed. Let%
\begin{equation}
w(x,t)\equiv a^{l}(x)u(x,t), \tag{3.3}%
\end{equation}
where $l$ is some positive constant to be determined later. Then the equation
in (1.1) is transferred into the following:%

\begin{align}
w_{t}-w_{xx}  &  =-2la^{-1}a_{x}w_{x}+\mu a^{(1-p)l}w^{p}-a^{1+l(1-q)}%
w^{q}\tag{3.4}\\
&  +[l^{2}a_{x}^{2}a^{-2}+la_{x}^{2}a^{-2}-la_{xx}a^{-1}]w.\nonumber
\end{align}
Given a number $\alpha\equiv(p-1)l,$then the following is obtained:%

\begin{align}
a^{\alpha}(w_{t}-w_{xx})  &  =-2la^{(p-1)l-1}a_{x}w_{x}+\mu w^{p}%
-u^{q-p-l^{-1}}w^{p+l^{-1}}\tag{3.5}\\
&  +a^{(p-1)l-1}[l^{2}a_{x}^{2}a^{-1}+la_{x}^{2}a^{-1}-la_{xx}]w.\nonumber
\end{align}
If a sufficiently large $l$ is chosen, then the coefficients of $w_{x},$
$w^{p+l^{-1}},$ and $w$ are smooth. Let the coefficients of $w$ be bounded by
$M$. Since $u$ blows up in a finite time, $u\neq$constant. Meanwhile, since
$u_{0}(x)\geqq0$, then $u(x,t)>0$ for ($x,t)\in\overline{D}\times(0,T^{\ast})$
by the weak and Hopf maximum principles. Without a loss of generality,
$u(x,t)>0$ is assumed for ($x,t)\in\overline{D}\times\lbrack0,T^{\ast}).$
Define $\beta\equiv\min_{(x,t)\in\overline{D}\times\lbrack0,T^{\ast}%
]}u(x,t)>0.$ Let $w^{+}$ be the largest root of the following equation:%

\[
\mu w^{p-1}-\beta^{q-p-l^{-1}}w^{p+l^{-1}-1}+M=0
\]
and $w^{\ast}$ be the following number:%

\[
w^{\ast}\equiv\max\{w^{+},\max_{x\in\overline{D}}w(x,0)\}>0.
\]
In the following, it will be shown that
\begin{equation}
w(x,t)\leqq w^{\ast},\text{ for (}x_{,}t)\in\overline{D}\times\lbrack
0,T^{\ast}) \tag{3.6}%
\end{equation}
and then the conclusion of the theorem can be implied. \ \ \ \ \ \ \ \ \ \ \ \ \ \ 

If (3.6) is not true, then there exists ($x_{0,}t_{0})\in\overline{D}%
\times(0,T^{\ast})$ such that $w_{0\equiv}w(x_{0},t_{0})>w^{\ast}$ $>0.$and
$w_{0}=\max_{\text{(}x_{,}t)\in\overline{D}\times\lbrack0,t_{0}]}.$ Therefore%

\begin{equation}
\mu w_{0}^{p-1}-\beta^{q-p-l^{-1}}w_{0}^{p+l^{-1}-1}+M<0. \tag{3.7}%
\end{equation}

If $x_{0}\in D,$ then we have%

\[
a^{\alpha}(x_{0})[w_{t}(x_{0},t_{0})-w_{xx}(x_{0},t_{0})]\geqq0\text{ and
}w_{x}(x_{0},t_{0})=0.
\]
Therefore from (3.5) we have%
\[
\mu w_{0}^{p-1}-u^{q-p-l^{-1}}w_{0}^{p+l^{-1}-1}+a^{(p-1)l-1}[l^{2}a_{x}%
^{2}a^{-1}+la_{x}^{2}a^{-1}-la_{xx}]\geqq0.
\]
However, this contradicts (3.7) by the definitions of $\beta,$ $M.$

If $x_{0}\in\partial D,$ then $a(x_{0})>0$ and so $\frac{\partial a(x_{0}%
)}{\partial n}\leqq0.$ Using (3.7) and the continuity of $w(x,t)$ and $a(x),$
there exists a $\delta$ $>0$ small enough such that%

\begin{align}
h(x,t)  &  \equiv\mu w^{p-1}-u^{q-p-l^{-1}}w^{p+l^{-1}-1}+a^{(p-1)l-1}%
[l^{2}a_{x}^{2}a^{-1}+la_{x}^{2}a^{-1}-la_{xx}]\tag{3.8}\\
&  \leqq\mu w^{p-1}-\beta^{q-p-l^{-1}}w^{p+l^{-1}-1}+M<0\text{ on }%
\overline{B(x_{0},\delta)\cap D}\times\lbrack t_{0}-\delta,t_{0}].\nonumber
\end{align}
and%
\begin{equation}
a(x)>0\text{ on on }\overline{B(x_{0},\delta)\cap D}. \tag{3.9}%
\end{equation}
Because $w\neq$constant and (3.8), (3.9), $\frac{\partial w(x_{0},t_{0}%
)}{\partial n}>0$ is obtained. However, the following conditions: $\left.
\frac{\partial u}{\partial n}\right|  _{\partial D}=0$ and $\frac{\partial
a(x_{0})}{\partial n}\leqq0,$ lead to a contradiction.

Second, problem ($\mathcal{III}$) is considered. Let $l\equiv(q-p)^{-1}$ and
then (3.5) becomes the following:%
\begin{align}
a^{\alpha}(w_{t}-w_{xx})  &  =-2la^{(p-1)l-1}a_{x}w_{x}+\mu w^{p}-w^{p+l^{-1}%
}\tag{3.10}\\
&  +a^{(p-1)l-1}[l^{2}a_{x}^{2}a^{-1}+la_{x}^{2}a^{-1}-la_{xx}]w.\nonumber
\end{align}
Under the condition: $\sigma\geqq2(q-p)(p-1)^{-1},$ the coefficients of
$\ w_{x},$ $w$ in (3.10) are smooth. The result is obtained by proceeding as
in problem ($\mathcal{II}$) and noting that $x_{0}$ above does not occur at
$\partial D$.
$\ \ \ \ \ \ \ \ \ \ \ \ \ \ \ \ \ \ \ \ \ \ \ \ \ \ \ \ \ \ \ \ \ \ \ \ \ \ \ \ \ \ \ \ \ \blacksquare
$

\bigskip

\bigskip

\section{Application: \ Diffusion-Induced Blowup}

This section attempts to demonstrate that the reaction diffusion system (1.6)
has the \textit{diffusion-induced blowup} phenomenon.

\bigskip

\bigskip

\textbf{Proposition 4.1. }\textit{Assume that q}$\geqq p>1$ $and$ $\sigma
\geqq2(q-p)(p-1)^{-1}.$ \textit{When} $n\geqq3,$ \textit{it is further assumed
that} $1<q<\frac{n+2}{n-2}.$ $\mathit{If}$ $\mu$ $\mathit{and}$ $u_{0}(x)$
\textit{are large enough}$,$\textit{ then} \textit{the solutions }$u(x,t)$
\textit{for the} \textit{Cauchy, Neumann, and Dirichlet problems of (1.6) blow
up in finite time.}

\bigskip

\bigskip

\textit{Proof. }\textrm{This follows from Theorem 1.2.
\ \ \ \ \ \ \ \ \ \ \ \ \ \ \ \ \ \ \ \ \ \ \ \ \ \ \ \ \ \ \ \ \ \ \ \ \ \ \ \ \ \ \ \ \ \ \ \ \ \ \ \ \ \ }%
\textbf{$\blacksquare$}

\bigskip

\bigskip

\textbf{Proposition 4.2. }\textit{Assume that f(x}$_{0,}$\textit{v) is a
smooth function such that the solution v(t) for (1.7) exists globally for any
}$\xi\geqq0$ \textit{(e.g., f(x}$_{0},$\textit{0)}$\geqq0$ $,$ \textit{f}%
$(x_{0,}v)\leqq0$ \textit{when v is large, or f(x}$_{0},v)$ \textit{is linear
in v), q}$>p,$\textit{\ q}$>$\textit{1, }$and$ $\sigma\geqq0$\textit{. If
f(x}$_{0,}m)\neq0,$ \textit{then the solution u for the kinetic system (1.7)
is global, \ bounded, \ and\ nonnegative for any nonnegative initial data
}$\eta$\textit{.}

\bigskip

\bigskip

\textit{Proof.} The argument in [\textbf{18}] is modified to obtain the global
existence and boundedness of $u$. The solutions \textit{v,u }\textrm{of (1.7)
are denoted by }\textit{v(t;}$\xi)$ and $u(t;\eta)$. If there exist $\xi
_{0},\eta_{0}$ such that $u(t;\eta_{0})$ blows up at finite time $T^{\ast}$,
from the phase plane analysis,\ then $v(T^{\ast};\xi_{0})=m.$ Notably, the
$v-$equation in (1.7) is autonomous. Thus, only three type of orbits of $v$
can occur: equilibrium, closed orbit, and strictly monotone orbit. Since $v$
is scalar, the $\omega-$limit set of $v$ consists of just one equilibrium
point; thus the orbit can not be the closed one. By $f(x_{0},m)\neq0,$ $v$ is
strictly monotone. Hence, $T^{\ast}$ is the only finite time such that $u$
blows up. By the comparison principle of the ordinary differential equation,
$u(t;\eta)$ also blows up at time $T^{\ast}$ for $\eta\geqq\eta_{0}$. Now, we
want to demonstrate the following fact: $u(t;\eta)$ tends to infinity
uniformly on $[0$,$T^{\ast}],$ as $\eta$ tends to infinity. Contrarily, assume
that there exists a large number $M>\eta_{0}$ satisfying: for any $\eta>M,$ a
finite time $t(\eta)>0$ exists such that $u(t(\eta);\eta)<M.$ Let $t_{M}%
(\eta)\equiv\min\{t:$ $u(t;\eta)=M\}$ for any $\eta>M.$ Because $u(t;\eta)$
blows up at $T^{\ast}$, by means of the comparison principle of ordinary
differential equation and sequentially compact property, an increasing
sequence \{$t_{M}(\eta_{i})\}_{i}$ which tends to some $t_{M}\in(0,T^{\ast})$
exists. The solution $u$ which passes through the point $(M,t_{M})$ in the
$u-t$ plane \ will blow up at some finite time less than $T^{\ast}.$ This
contradicts the fact: ``$T^{\ast}$ is the only finite time such that $u$ blows up.''

Set $w:$ $=u^{-(q-1)}.$ The following equation for \textit{w }is then
obtained:
\begin{equation}
w_{t}=-(q-1)\mu w^{\frac{q-p}{q-1}}+(q-1)\left|  m-v(t;\xi_{0})\right|
^{\sigma}. \tag{4.1}%
\end{equation}

From above, a family of $w(t;\eta^{-(q-1)})$ is obtained, which tends to zero
uniformly on [$0$,$T^{\ast}$], as $\eta$ tends to infinity.\textit{\ }%
\textrm{Since }$\mathit{f(x}_{0},m)\neq0,$ some $\widetilde{t}$ $\in
\lbrack0,T$* ) exists, such that $v(\widetilde{t};\xi_{0})\neq m$ and
then\ $w_{t}>0$ at the point ($\widetilde{t,}0)$ in the t-w plane by (4.1).
This will contradict the continuous dependence property of $w(t;\eta
^{-(q-1)})$ near the point ($\widetilde{t,}0)$ in the t-w plane. Therefore
\textit{v} exists globally.

From \textrm{\ }$\mathit{f(x}_{0},m)\neq0,$ it can be easily obtained that
\textit{v(t}$_{k}$\textit{;}$\xi)$ $\nrightarrow m$ \ for any sequence
$t_{k}\rightarrow$ $\infty.$ Thus, $\mathit{u}$\textit{ }is\textit{ }also bounded.

The proof is
finished.\ \ \ \ \ \ \ \ \ \ \ \ \ \ \ \ \ \ \ \ \ \ \ \ \ \ \ \ \ \ \ \ \ \ \ \ \ \ \ \ \ \ \ \ \ \ \ \ \ \ \ \ \ \ \ \ \ \ \ \ \ \ \ \ \ \ \ \ \ \ \ \ \ \textbf{$\blacksquare
$}

\bigskip

\bigskip

Some examples illustrating the \textit{diffusion-induced blowup} phenomenon
are as follows:

\textbf{( i )}.\textit{\ f(x,v)} : =\textit{\ }$\lambda$\textit{v,} \ where
$\lambda>0$ is an eigenvalue of -$\Delta$ for the homogeneous Dirichlet, or
Neumann problem in $D.$ Clearly, the solution\textit{\ v(t)} of (1.7) exists
globally but is not bounded. In addition, based on the eigenfunction's
property, it is possible to find \textit{x}$_{0}$$\in$ $D$ such that
\textit{f(x}$_{0},v_{0}($\textit{x}$_{0}$\textit{)) =}$\lambda v_{0}%
(\mathit{x}_{0}\mathit{)}\neq0,$ \ where \textit{v}$_{0}$ is the eigenfunction
with respect to $\lambda.$

\textbf{( ii )}. \textit{\ f(x,v) }: = $\lambda v-h(x)v^{l}$, \ where
$\lambda>0$, $l>1$ are constants, and \textit{h, } not identical to zero, is a
nonnegative smooth function with isolated zero points. From [\textbf{22}%
,\textbf{Theorem 2}] and maximum principle, an unique nonnegative solution
\textit{v}$_{0}$ of the elliptic Dirichlet problem and \textit{x}$_{0}\in D$
can be found such that \textit{f(x}$_{0},$\textit{\ v}$_{0}(x_{0}))\neq0,$
\ for $\lambda>\lambda_{0},$ \ where $\lambda_{0}$ is the first eigenvalue of
\ -$\Delta$ for the homogeneous Dirichlet problem in $D$. Clearly, the
solution $v(t)$ of (1.7) exists globally, and is also bounded with the further
restriction that ``$h(x)$ is positive.''

\textbf{( iii )}. Take the same\textit{\ f(x,v)} as in ( \textbf{ii }), but
where ``\textit{h} is not a constant.'' Also from [\textbf{22},Theorem 3] and
maximum principle, a unique nonnegative solution \textit{v}$_{0}$ of the
elliptic Neumann problem and \textit{x}$_{0}\in D$ can be found such that
\textit{f(x}$_{0},$\textit{\ v}$_{0}(x_{0}))\neq0,$ \ for $\lambda>0$. The
conclusion for $v(t)$ is the same as in ( ii ).

\bigskip

\bigskip

Therefore, under the conditions ``\textit{\ }$q>p>1,1<q<\frac{n+2}{n-2}$
$(n\geqq3)$ $or$ $n=1,2,$ $and$ $\sigma\geqq2(q-p)(p-1)^{-1}",$ the above
examples have the so-called \textit{diffusion-induced blowup} phenomenon,
particularly in the case\textit{\ f(x,v) = f(v).}

\bigskip

\bigskip

\textbf{Remark. }\textrm{When we }take \textit{d}$_{1}=d_{2}$ in the examples
herein, the question in [\textbf{20}] as mentioned in the Introduction has
been resolved.

\ \ \ \ \ \ \ \ \ \ \ \ \ \ \ \ \ \ \ \ \ \ \ \ \ \ \ \ \ \ \ \ \ \ \ \ \ \ \ \ \ \ \ \ \ \
\[
\mathbf{\ Acknowledgements}%
\]

The author want to thank Prof. Shao-Shiung Lin for his kindly, useful suggestions.

\end{document}